\documentclass[11pt,a4paper]{article}

\usepackage[T1]{fontenc}
\usepackage[utf8]{inputenc}
\usepackage{lmodern}
\usepackage{microtype}
\usepackage{amsmath,amssymb,amsthm,mathtools,amscd}
\usepackage{enumitem}
\usepackage{hyperref}
\usepackage[nameinlink,noabbrev]{cleveref}
\usepackage[margin=30mm]{geometry}

\hypersetup{
  colorlinks=true,
  linkcolor=blue,
  citecolor=blue,
  urlcolor=blue,
  pdftitle={The Seifert--van Kampen Theorem: A Self-Contained Exposition},
  pdfauthor={}
}

\theoremstyle{plain}
\newtheorem{theorem}{Theorem}[section]

\theoremstyle{definition}
\newtheorem{definition}[theorem]{Definition}

\theoremstyle{remark}
\newtheorem{remark}[theorem]{Remark}

\newcommand{\pione}{\pi_1}

\newcommand{\Ker}{\mathrm{Ker}}

\title{\textbf{The Seifert--van Kampen Theorem}\\
\large A self-contained exposition with applications}
\author{Djamel AIT-AKLI}
\date{\today}

\begin{document}
\maketitle

\begin{abstract}
The Seifert--van Kampen theorem computes the fundamental group of a
space from an open cover by two path-connected subspaces, provided
their intersection is path-connected.  This article gives a
self-contained account of the theorem, explains its universal-property
formulation as a pushout in the category of groups, and derives several
standard computations.  Particular attention is paid to base points,
the role of path-connectedness, and the passage from a topological
decomposition to a group presentation.  The purpose is expository:
the results are classical, and no claim of a new theorem is made.
\end{abstract}

\noindent\textbf{2020 Mathematics Subject Classification:}
55Q05, 20E05.\\
\textbf{Keywords:} fundamental group, free product with amalgamation,
Seifert--van Kampen theorem, wedge sum, cell attachment.

\tableofcontents

\section{Introduction}

The fundamental group is one of the first invariants that detects
nontrivial global topology.  Its definition involves loops based at a
point and homotopies between such loops.  Directly classifying all
loops in a space is often difficult.  The Seifert--van Kampen theorem
provides a powerful alternative: if a space is decomposed into two
suitable open pieces, then its fundamental group can be reconstructed
from the fundamental groups of the pieces and their intersection.

The algebraic construction appearing in the theorem is the
\emph{free product with amalgamation}.  If two groups receive
homomorphisms from a third group, their pushout is the group obtained
by imposing precisely the relations that identify the two images of
each element of the common group.  The theorem says that this
algebraic pushout models the fundamental group of the union.

We work throughout with based spaces and continuous maps.  Unless
otherwise stated, all spaces are assumed to be locally path-connected
where this is needed to ensure that the chosen open sets and their
intersection have the stated path-connectedness properties.  The
precise hypotheses used in the main theorem are given explicitly.

\section{Fundamental groups and group pushouts}

\begin{definition}
Let $(X,x_0)$ be a pointed topological space.  A loop at $x_0$ is a
continuous map $\gamma:[0,1]\to X$ such that
$\gamma(0)=\gamma(1)=x_0$.  Two loops are based-homotopic if they are
homotopic through loops that keep the endpoints equal to $x_0$.
The set of based-homotopy classes of loops, with multiplication
induced by concatenation, is the fundamental group
$\pione(X,x_0)$.
\end{definition}

For a continuous pointed map $f:(X,x_0)\to(Y,y_0)$, composition
induces a homomorphism
\[
 f_*:\pione(X,x_0)\longrightarrow\pione(Y,y_0).
\]
If $A\subseteq X$ and $x_0\in A$, the inclusion induces a homomorphism
$\pione(A,x_0)\to\pione(X,x_0)$.

\begin{definition}
Let $G_0,G_1,G_2$ be groups and let
$\varphi_1:G_0\to G_1$ and $\varphi_2:G_0\to G_2$ be homomorphisms.
A pushout of this diagram in the category of groups is a group $P$
together with homomorphisms $j_1:G_1\to P$ and $j_2:G_2\to P$
such that $j_1\varphi_1=j_2\varphi_2$, and with the following
universal property: for every group $H$ and homomorphisms
$u_1:G_1\to H$, $u_2:G_2\to H$ satisfying
$u_1\varphi_1=u_2\varphi_2$, there exists a unique homomorphism
$u:P\to H$ such that $u j_i=u_i$ for $i=1,2$.
\end{definition}

The pushout exists and is denoted $G_1 *_{{G_0}} G_2$.  Concretely,
it is the quotient of the free product $G_1*G_2$ by the normal
subgroup generated by
\[
 \varphi_1(g)\varphi_2(g)^{-1},\qquad g\in G_0,
\]
where the two factors are understood as their canonical images in
the free product.  Thus
\[
G_1 *_{{G_0}} G_2
\cong
(G_1*G_2)\big/
\left\langle\!\left\langle
\varphi_1(g)\varphi_2(g)^{-1}:g\in G_0
\right\rangle\!\right\rangle.
\]

\section{Statement of the theorem}

\begin{theorem}[Seifert--van Kampen]
\label{thm:svk}
Let $X=U\cup V$, where $U$ and $V$ are open and path-connected
subspaces of $X$.  Suppose that $U\cap V$ is nonempty and
path-connected.  Choose a base point $x_0\in U\cap V$.
Then the homomorphisms induced by inclusion,
\[
\pione(U\cap V,x_0)\longrightarrow\pione(U,x_0),
\qquad
\pione(U\cap V,x_0)\longrightarrow\pione(V,x_0),
\]
make the following diagram a pushout in the category of groups:
\[
\begin{CD}
\pione(U\cap V,x_0) @>>> \pione(U,x_0)\\
@VVV @VVV\\
\pione(V,x_0) @>>> \pione(X,x_0).
\end{CD}
\]
Equivalently,
\[
\pione(X,x_0)\cong
\pione(U,x_0)*_{\pione(U\cap V,x_0)}\pione(V,x_0).
\]
\end{theorem}

The openness assumptions are important in the standard proof because
they allow a path or homotopy on a compact interval or square to be
subdivided into sufficiently small pieces lying in one member of the
cover.  There are more general versions, but the formulation above is
the one used in the applications below.

\section{Proof}

\begin{proof}
We prove the universal property.  Let
\[
G=\pione(U,x_0)*_{\pione(U\cap V,x_0)}\pione(V,x_0).
\]
The inclusions of $U$ and $V$ into $X$ induce homomorphisms to
$\pione(X,x_0)$.  Their restrictions to $\pione(U\cap V,x_0)$ agree,
so the universal property of the pushout gives a homomorphism
\[
\Phi:G\longrightarrow\pione(X,x_0).
\]
We show that $\Phi$ is an isomorphism.

\subsection{Surjectivity}

Let $[\gamma]\in\pione(X,x_0)$ be represented by a loop
$\gamma:[0,1]\to X$.  Since $[0,1]$ is compact and
$\{\gamma^{-1}(U),\gamma^{-1}(V)\}$ is an open cover, the Lebesgue
number lemma gives a subdivision
\[
0=t_0<t_1<\cdots<t_n=1
\]
such that each subinterval $[t_{i-1},t_i]$ is mapped entirely into
either $U$ or $V$.

At a subdivision point where the chosen member changes, the point
$\gamma(t_i)$ lies in $U\cap V$.  Because $U\cap V$ is path-connected,
choose a path in $U\cap V$ from $x_0$ to each such transition point.
Also use the constant path at $x_0$ at the two endpoints.  These
connecting paths express the original loop, up to based homotopy, as
a product of loops each of which lies entirely in $U$ or entirely in
$V$: one inserts a path followed by its reverse at each subdivision
point.  Consequently $[\gamma]$ is in the image of $\Phi$.
Thus $\Phi$ is surjective.

\subsection{Injectivity}

We use a standard subdivision argument for homotopies.  Suppose an
element of $G$ is represented by a product
\[
[g_1][h_1]\cdots[g_r][h_r],
\]
where each $g_i$ is a loop in $U$ and each $h_i$ is a loop in $V$,
with the usual interpretation that factors may be absent.  Suppose
its image under $\Phi$ is the identity.  Then the corresponding
concatenated loop $\gamma$ is based-null-homotopic in $X$; let
\[
H:[0,1]\times[0,1]\longrightarrow X
\]
be a based null-homotopy, with $H(s,0)=\gamma(s)$,
$H(s,1)=x_0$, and $H(0,t)=H(1,t)=x_0$.

The inverse images $H^{-1}(U)$ and $H^{-1}(V)$ form an open cover of
the compact square.  By the Lebesgue number lemma, choose a sufficiently
fine rectangular subdivision so that the image of every small
rectangle lies in $U$ or in $V$.  Along an edge shared by rectangles
of different types, its image lies in $U\cap V$.

Reading the boundary of each small rectangle gives a relation among
the edge paths in the corresponding group.  If the rectangle maps
into $U$, this is a relation in $\pione(U,x_0)$ after transporting
the boundary to the base point by paths in the connected intersection
where needed; if it maps into $V$, it is a relation in
$\pione(V,x_0)$.  The relations on common edges agree in the
amalgamated product because the two copies of every loop in
$U\cap V$ are identified there.  Multiplying the rectangle
relations over the grid cancels each interior edge with the same edge
traversed in the opposite direction.  The remaining outer boundary
is the original loop together with constant paths.  It follows that
the original word is trivial in $G$.

Therefore $\Ker(\Phi)$ is trivial, and $\Phi$ is injective.  This
proves \cref{thm:svk}.
\end{proof}

\begin{remark}
The injectivity argument is often presented using a grid subdivision
of the homotopy square.  To make every base-point transport explicit,
one may choose paths in $U\cap V$ from $x_0$ to the vertices lying on
the interfaces and conjugate the corresponding edge loops by these
paths.  The resulting relations are precisely the amalgamation
relations.  This is the geometric content of the universal-property
proof.
\end{remark}

\section{Applications}

\subsection{The fundamental group of a wedge of two circles}

Let $X=S^1\vee S^1$ be the wedge of two circles, with wedge point
$x_0$.  Choose open neighborhoods $U$ and $V$ of the two circles,
slightly enlarged around the wedge point, so that $U$ deformation
retracts onto the first circle, $V$ deformation retracts onto the
second circle, and $U\cap V$ is contractible.  The theorem gives
\[
\pione(S^1\vee S^1,x_0)
\cong \pione(S^1,x_0)*\pione(S^1,x_0)
\cong \mathbb Z*\mathbb Z.
\]
The free product of two infinite cyclic groups is the free group
$F_2$ on two generators.  More generally, the wedge of $n$ circles
has fundamental group the free group $F_n$.

\subsection{The fundamental group of a circle}

For completeness, the standard computation
$\pione(S^1,1)\cong\mathbb Z$ can be obtained from the covering map
$p:\mathbb R\to S^1$, $p(t)=e^{2\pi i t}$.  A loop based at $1$
has a unique lift starting at $0$, and its endpoint is an integer.
The endpoint is invariant under based homotopy and is additive under
concatenation.  Every integer occurs, and a loop whose lift has
endpoint $0$ is null-homotopic.  Hence the endpoint map is an
isomorphism with $\mathbb Z$.

\subsection{The fundamental group of a sphere}

For $n\geq 2$, write $S^n=U\cup V$, where $U$ and $V$ are the
complements of two distinct points, slightly adjusted as open
hemispherical neighborhoods.  Each of $U$ and $V$ is contractible,
and their intersection is path-connected for $n\geq 2$ (it deformation
retracts onto $S^{n-1}$).  Choosing $x_0\in U\cap V$, the theorem yields
\[
\pione(S^n,x_0)\cong 1*_{\pione(S^{n-1},x_0)}1=1,
\qquad n\geq 2.
\]
Thus every sphere of dimension at least two is simply connected.
For $n=2$, the intersection is homotopy equivalent to $S^1$, but
both maps to the trivial groups are trivial, so the pushout is still
trivial.

\subsection{A presentation for a surface group}

The closed orientable surface of genus $g$ has a CW decomposition
with one vertex, $2g$ one-cells
$a_1,b_1,\ldots,a_g,b_g$, and one two-cell attached along the word
\[
[a_1,b_1]\cdots[a_g,b_g],
\qquad [a,b]=aba^{-1}b^{-1}.
\]
Applying the van Kampen theorem to the one-skeleton and the attached
two-cell gives
\[
\pione(\Sigma_g)\cong
\left\langle a_1,b_1,\ldots,a_g,b_g
\ \middle|\
[a_1,b_1]\cdots[a_g,b_g]=1
\right\rangle.
\]
This illustrates a central use of the theorem: attaching a
two-dimensional cell adds a relation represented by its attaching
loop.

\section{Conceptual interpretation}

The theorem is a compatibility statement between topology and
category theory.  The union $X=U\cup V$ is a colimit-like construction
in topology, while the fundamental-group functor converts the cover
into a pushout diagram of groups under the stated hypotheses.

It is important not to interpret the result as saying that the
fundamental group preserves every union.  The hypotheses on the
intersection, the base point, and the open cover are part of the
theorem.  In particular, if the intersection is not path-connected,
one must use a more elaborate formulation, often involving groupoids.
The fundamental groupoid version is more flexible because it can
retain several base points simultaneously.

\section{Further directions}

Several natural extensions place the theorem in a broader homotopical
framework.

\begin{enumerate}[label=(\arabic*)]
\item \textbf{The groupoid van Kampen theorem.}
  Replacing one base point by a set of base points gives a theorem
  that handles disconnected intersections and is often more natural
  for complicated decompositions.

\item \textbf{CW complexes and presentations.}
  Repeated application of van Kampen gives presentations of fundamental
  groups from the one- and two-skeleta of a CW complex.  Cells of
  dimension at least three do not change the fundamental group.

\item \textbf{Homotopy colimits.}
  The pushout formula is the degree-one shadow of broader
  homotopy-colimit principles.  In modern homotopy theory, one studies
  conditions under which homotopy-invariant functors turn suitable
  homotopy pushouts into homotopy pullbacks or other computable
  constructions.

\item \textbf{Higher homotopy.}
  The fundamental group is only the first nontrivial homotopy
  invariant.  Mayer--Vietoris sequences and homotopy excision provide
  complementary tools for understanding higher-dimensional
  information.
\end{enumerate}

\section{Conclusion}

The Seifert--van Kampen theorem transforms a topological decomposition
into an algebraic pushout.  Its practical strength comes from the
combination of two ideas: loops can be subdivided so that each piece
lies in one member of an open cover, and the ambiguity at transitions
is controlled by the fundamental group of the intersection.  This
makes the theorem indispensable for computing fundamental groups of
graphs, surfaces, and spaces built by attaching cells.


\begin{thebibliography}{9}

\bibitem{Hatcher}
A. Hatcher,
\emph{Algebraic Topology},
Cambridge University Press, Cambridge, 2002.
Available online at
\url{https://pi.math.cornell.edu/~hatcher/AT/AT.pdf}.

\bibitem{May}
J. P. May,
\emph{A Concise Course in Algebraic Topology},
University of Chicago Press, Chicago, 1999.

\bibitem{Massey}
W. S. Massey,
\emph{A Basic Course in Algebraic Topology},
Graduate Texts in Mathematics, vol. 127,
Springer-Verlag, New York, 1991.

\bibitem{tomDieck}
T. tom Dieck,
\emph{Algebraic Topology},
European Mathematical Society, Z\"urich, 2008.

\end{thebibliography}
\end{document}